\documentstyle[amscd]{amsart}
\numberwithin{equation}{section}
\theoremstyle{plain}
\newtheorem{thm}{Theorem}[section]
\newtheorem{cor}[thm]{Corollary}
\newtheorem{lem}[thm]{Lemma}
\newtheorem{prop}[thm]{Proposition}
\begin{document}
\title{
A class of simple $C^*$-algebras arising from certain nonsofic subshifts
}
\author{Kengo Matsumoto}
\address{ 
Department of Mathematical Sciences, 
Yokohama City  University,
22-2 Seto, Kanazawa-ku, Yokohama, 236-0027 Japan}
\email{kengo@@yokohama-cu.ac.jp}
\maketitle
\begin{abstract}
We present a class of subshifts $Z_N, N = 1,2,\dots$ whose associated 
$C^*$-algebras ${\cal O}_{Z_N}$ are simple, purely infinite 
and not stably isomorphic to any Cuntz-Krieger algebra nor to Cuntz algebra.
The class of the subshifts is the first examples 
whose associated $C^*$-algebras are 
not stably isomorphic to any Cuntz-Krieger algebra nor to Cuntz algebra.
The subshifts $Z_N$ are coded systems whose languages are context free.
We compute the topological entropy for the subshifts and
 show that
KMS-state for gauge action on the associated $C^*$-algebra ${\cal O}_{Z_N}$
exists
 if and only if the logarithm of the inverse temperature is
  the topological entropy for the subshift $Z_N$,
  and the corresponding KMS-state is unique.
 \end{abstract}

\keywords{
Simple $C^*$-algebra, subshift, entropy, KMS-state.
}

\subjclass{ 
Primary 46L35, 37B10; Secondary 54H20,54C70
}

\def\Zp{{ {\Bbb Z}_+ }}
\def\OL{{{\cal O}_{{\frak L}}}}
\def\OZN{{{\cal O}_{Z_N}}}
\def\M{{{\cal M}}}
\def\N{{{\cal M}}}
\def\H{{{\cal H}}}
\def\K{{{\cal K}}}
\def\P{{{\cal P}}}
\def\Q{{{\cal Q}}}
\def\A{{{\cal A}}}
\def\B{{{\cal B}}}
\def\R{{{\cal R}}}
\def\S{{{\cal S}}}
\def\sms{{{symbolic  matrix system }}}
\def\smss{{{symbolic  matrix systems }}}
\def\nnms{{{nonnegative matrix system }}}
\def\nnmss{{{nonnegative  matrix systems }}}
\def\Ext{{{\operatorname{Ext}}}}
\def\Aut{{{\operatorname{Aut}}}}
\def\Hom{{{\operatorname{Hom}}}}
\def\Ker{{{\operatorname{Ker}}}}
\def\id{{{\operatorname{id}}}}
\def\FKL{{ {\cal F}_k^{l} }}
\def\FKI{{{\cal F}_k^{\infty}}}
\def\FLI{{ {\cal F}_{\Lambda} }}
\def\SES{{S_{\mu}E_l^iS_{\mu}^*}}


\section{Introduction}

Let $\Sigma$ be a finite set with its discrete topology, 
that is called an alphabet.
Each element of $\Sigma$ is called a symbol.
Let $\Sigma^{\Bbb Z}$ 
be the infinite product space 
$\prod_{i\in {\Bbb Z}}\Sigma_{i},$ 
where 
$\Sigma_{i} = \Sigma$,
 endowed with the product topology.
 The transformation $\sigma$ on $\Sigma^{\Bbb Z}$ 
given by 
$
\sigma( (x_i)_{i \in {\Bbb Z}})
 = (x_{i+1})_{i \in {\Bbb Z}}
$ 
is called the full shift over $\Sigma$.
 Let $\Lambda$ be a shift invariant closed subset of $\Sigma^{\Bbb Z}$ i.e. 
 $\sigma(\Lambda) = \Lambda$.
  The topological dynamical system 
  $(\Lambda, \sigma\vert_{\Lambda})$
   is called a subshift or a symbolic dynamical system.
   It is written as $\Lambda$ for brevity.
  Theory of symbolic dynamical system gives a basic method to study general dynamical systems (cf.\cite{LM}).
     It also has significant uses in coding for information theory.  
   There is a class of subshifts called sofic shifts,
that contains the topological Markov shifts.
Sofic shifts are presented by finite square matrices with entries 
in formal sums of symbols.
Such a matrix is called a symbolic matrix.
It is an equivalent object to a finite labeled graph
 called a $\lambda$-graph.  
In \cite{Ma5}, 
the  author  
has introduced the notions of 
symbolic matrix system and
$\lambda$-graph system.
They are presentations of subshifts
and generalizations  of 
 symbolic matrices and
$\lambda$-graphs respectively.
A symbolic matrix  system $(\M,I)$ consists of
a sequence  of pairs $(\M_{l,l+1}, I_{l,l+1}), l \in \Zp$
of rectangular symbolic matrices
$\M_{l,l+1}$
and
rectangular $\{0,1\}$-matrices
$I_{l,l+1}$, 
where $\Zp$ denotes the set of all 
nonnegative integers.
Both the matrices $\M_{l,l+1}$ and $I_{l,l+1}$ have the same size
for each $l \in \Zp.$
The column size of $\M_{l,l+1}$ is the same as the row size of 
$
\M_{l+1,l+2}.
$
They satisfy the following commutation relations as symbolic matrices
\begin{equation}
I_{l,l+1} \M_{l+1,l+2} = \M_{l,l+1}I_{l+1,l+2}, 
\qquad
l \in \Zp. 
\end{equation}
We further assume  that
for $i$ 
there exists $j$ such that
the $(i,j)$-component $I_{l,l+1}(i,j) =1,$
and for  $j$ 
there uniquely exists $i$ such that
$I_{l,l+1}(i,j) =1$.
A $\lambda$-graph system 
$ {\frak L} = (V,E,\lambda,\iota)$
consists of a vertex set 
$V = V_0\cup V_1\cup V_2\cup\cdots$, an edge set 
$E = E_{0,1}\cup E_{1,2}\cup E_{2,3}\cup\cdots$, 
a labeling map
$\lambda: E \rightarrow \Sigma$
and a surjective map
$\iota_{l,l+1}: V_{l+1} \rightarrow V_l$ for each
$l\in \Zp.$
An edge $e \in E_{l,l+1}$ has its source vertex $s(e)$ in $V_l,$
its terminal vertex $t(e)$ in $V_{l+1}$ 
and its label $\lambda(e)$ in $\Sigma$.
For a symbolic matrix system  $(\M,I)$,
the labeled edges from a vertex $v_i^l \in V_l$ 
to a vertex  $v_j^{l+1}\in V_{l+1}$ are given by the symbols appearing in the
$(i,j)$-component $\M_{l,l+1}(i,j)$
of $\M_{l,l+1}$.
The matrix $I_{l,l+1}$ 
defines a surjection $\iota_{l,l+1}$ 
from $V_{l+1}$ to $V_l$ for each $l \in \Zp.$
By this observation, 
the symbolic matrix systems and the $\lambda$-graph systems 
 are the same objects.
 They give rise to subshifts by taking the set of all label sequences
 appearing in the labeled Bratteli diagram.
For a symbolic matrix system
$(\M,I)$, 
let
$M_{l,l+1}$ be the nonnegative rectangular matrix obtained from 
$\M_{l,l+1}$ by setting all the symbols equal to $1$ 
for each $l \in \Zp$.
Then the resulting pair
$(M,I)$ satisfies the  following relations by (1.1)
\begin{equation}
I_{l,l+1} M_{l+1,l+2} = M_{l,l+1}I_{l+1,l+2}, 
\qquad
l \in \Zp. 
\end{equation}
We call
$(M,I)$
the nonnegative matrix system
for $(\M,I)$.
Let $m(l)$ be the row size of the matrix $I_{l,l+1}$ for each $l\in \Zp$.
Let 
${\Bbb Z}_{I^t}$
be the abelian group defined by the inductive limit
$
{\Bbb Z}_{I^t} = 
\varinjlim \{ I_{l,l+1}^t : {\Bbb Z}^{m(l)} \rightarrow {\Bbb Z}^{m(l+1)} \}.
$
The sequence $ M_{l,l+1}^t, l \in \Zp$ of the transposes of $M_{l,l+1}$
 naturally acts on
 ${\Bbb Z}_{I^t}$ by the relation (1.2),  
 that is denoted by $\lambda_{(M,I)}$.
 The K-groups for $(M,I)$  have been  defined in \cite{Ma5} as:
\begin{equation*}
K_0(M,I)  = {\Bbb Z}_{I^t} / (\id - \lambda_{(M,I)}){\Bbb Z}_{I^t},\qquad
K_1(M,I)   = \Ker(\id - \lambda_{(M,I)} ) 
              \text{ in } {\Bbb Z}_{I^t}.
\end{equation*}

In \cite{Ma9}, $C^*$-algebra ${\cal O}_{\frak L}$
associated with a $\lambda$-graph system ${\frak L}$ has been introduced.
The $C^*$-algebras are generalizations of the Cuntz-Krieger algebras and the $C^*$-algebras associated with subshifts
 (\cite{CK}, \cite{Ma}, cf. \cite{CaMa}).
They are universal unique concrete $C^*$-algebras generated 
by finite families of partial isometries and  sequences of projections subject to certain operator relations  encoded by structure of the $\lambda$-graph systems.
Let $(M,I)$ be the nonnegative matrix system for the symbolic matrix system 
$(\M,I)$ of ${\frak L}$. 
The K-theory formulae for ${\cal O}_{\frak L}$  
have been obtained  as in the following way:
\begin{equation}
K_0({\cal O}_{\frak L})  = K_0(M,I), \qquad
K_1({\cal O}_{\frak L}) = K_1(M,I).
\end{equation}
There is a canonical method to construct a $\lambda$-graph system 
${\frak L}$ from a given subshift $\Lambda$.
The $\lambda$-graph system is called the canonical $\lambda$-graph system
for $\Lambda$ and written as 
${\frak L}^\Lambda$.
The $C^*$-algebra 
${\cal O}_{{\frak L}^\Lambda}$
associated with the canonical $\lambda$-graph system
${\frak L}^\Lambda$ coincides with the $C^*$-algebra
${\cal O}_\Lambda$ associated with subshift $\Lambda$
( cf. \cite{CaMa}, \cite{Ma}).
If in particular $\Lambda$ is a topological Markov shift 
$\Lambda_A$ for a finite square matrix $A$ with entries in 
$\{0,1\}$, 
the $C^*$-algebra ${\cal O}_{{\frak L}^{\Lambda^A}}$
is canonically isomorphic to the Cuntz-Krieger algebra 
${\cal O}_A$.

In this paper,
we present a class of nonsofic subshifts 
$Z_N, N \in {\Bbb N}$ 
whose associated 
$C^*$-algebras ${\cal O}_{Z_N}$ are simple, purely infinite 
and not stably isomorphic to any Cuntz-Krieger algebra nor to Cuntz algebra.
The subshifts $Z_N, N \in {\Bbb N}$ 
are coded systems whose languages are context free (cf. \cite{BH}, \cite{ChS}, \cite{HU}).
In studying a topological dynamical system, 
the topological entropy is very important quantity to measure 
complexity for the topological dynamical system.
In \cite{EFW2}, 
it was shown that
the topological entropy for irreducible Markov shifts appear as 
the logarithm of the inverse temperature 
for admitting KMS-state for gauge action on the corresponding Cuntz-Krieger algebras.
This result has been generalized to more general subshifts in \cite{MWY}
(cf. \cite{LaNe}, \cite{PWY}).
Corresponding to these results, 
we will  
 compute the topological entropy for the subshifts $Z_N, N\in {\Bbb N}$ 
and
 show that
KMS-state for gauge action on the associated $C^*$-algebra ${\cal O}_{Z_N}$
exists
 if and only if the logarithm of the inverse temperature is
  the topological entropy for the subshift $Z_N$,
  and the corresponding KMS-state is unique.
We will prove 
\begin{thm}[Theorem 3.7, Theorem 4.8 and Theorem 5.15]
\qquad
\begin{enumerate}
\renewcommand{\labelenumi}{(\roman{enumi})}
\item
For $N \in {\Bbb N}$, 
the $C^*$-algebra $\OZN$ 
associated with the subshift 
$Z_N$ is simple and purely infinite.
It is the universal concrete 
$C^*$-algebra generated by two isometries $T_1, T_2$
and $N$ partial isometries $S_i, i=1,\dots, N$
satisfying the relations:
\begin{align}
T_1 T_1^* & +  T_2 T_2^* + \sum_{j=1}^N S_j S_j^* =1,\\
S_i^*S_i  =& 1 - \sum_{k=1}^{\infty} \sum_{m=0, k \ne 2m}^{k} \sum_{j=1}^N
           T_{2^m1^{k-m}} S_jS_j^* T_{2^m1^{k-m}}^*, \qquad i=1,\dots,N 
\end{align}
where
 $ T_{2^m1^{k-m}}$ 
 denotes
  $
  \underbrace{T_2 \cdots T_2}_{\text{$m$ times}}
  \underbrace{T_1 \cdots T_1}_{\text{$k-m$ times}}.  
 $
The infinite sum of the right hand side of the relation (1.5) is taken
under strong operator topology on a Hilbert space.
\item
The K-groups for the $C^*$-algebra $\OZN$ 
are as follows:
$$
K_0(\OZN) = {\Bbb Z}/N{\Bbb Z} \oplus {\Bbb Z}
\qquad
\text{and} 
\qquad
K_1(\OZN) = 0.
$$
The position of the unit $[1]$ in 
$K_0(\OZN) = {\Bbb Z}/N{\Bbb Z} \oplus{\Bbb Z} $
is 
$0$.
Hence $\OZN$ is not  stably isomorphic to any  Cuntz-Krieger algebra
nor  to the Cuntz algebra ${\cal O}_{\infty}$ for $N \ge 2$.
\item
There is a
KMS-state for gauge action on $\OZN$ 
if and only if the logarithm of the inverse temperature is $\log \beta_N$: 
the topological entropy of the subshift $Z_N$,
where $\beta_N$ is the unique solution $\beta$ satisfying $\beta > N$ 
of the equation
\begin{equation}
 \beta^6 -(N+3)\beta^5 +(3N+1)\beta^4- 2(N-1)\beta^3 -(N+2)\beta^2 + N-1 =0.
\end{equation}
The admitted KMS-state is unique.
\end{enumerate}
\end{thm}

The value $\beta_N$ is increasing on $N$ and satisfies
$N < \beta_N < N+1$ such that
\begin{equation*}
\lim_{N \to \infty} \frac{\beta_N}{N} =1.
\end{equation*}
The class of  the subshifts $Z_N, N\in {\Bbb N}$
is the first examples 
whose associated $C^*$-algebras are 
not stably isomorphic to any Cuntz-Krieger algebra nor to Cuntz algebra.
For $N = 1$, the subshift $Z_1$ is nothing but the subshift $Z$ 
named as the context free shift in \cite[Example 1.2.9]{LM}
and the associated $C^*$-algebra ${\cal O}_{Z_1}$ is the $C^*$-algebra ${\cal O}_Z$
studied in \cite{Ma3}.
For other type of coded systems, see for example \cite{KM4}.

We will finally mention  an application of our discussions to 
a classification problem in the theory of symbolic dynamical systems.
By \cite{Ma7}, \cite{Ma8},
the K-groups $K_i({\cal O}_\Lambda)$ 
are invariant under not only topological conjugacy but also
flow equivalence of subshifts $\Lambda$.
Hence we know that the subshifts
$Z_N, N\in {\Bbb N}$ are not flow equivalent to each other.

\section{Subshifts and the $C^*$-algebras associated with $\lambda$-graph systems}

For a subshift $\Lambda$ over alphabet 
$\Sigma$,
we denote by 
$X_{\Lambda} \subset \Sigma^{\Bbb N}$ 
the set of all right-infinite sequences that appear in $\Lambda$.
The topological dynamical system 
$(X_{\Lambda},\sigma)$ 
is called 
the one-sided subshift for $\Lambda$. 
   A finite sequence $\mu = (\mu_1,...,\mu_k) $ of elements $\mu_j \in \Sigma$
 is called a word. 
  We denote by $|\mu|$ the length $k$ of $\mu$. 
  A block $\mu = (\mu_1,...,\mu_k)$ is said to  appear in 
$x=(x_i) \in \Sigma^{\Bbb Z}$ if $x_m = \mu_1,..., x_{m+k-1} = \mu_k$  
for some $m \in {\Bbb Z}$.
For  a number $k\in {\Bbb N},$
let
$
B_k(\Lambda)
$
be the set of all words of length  
$k$
 in 
$ \Sigma^{\Bbb Z}$
appearing in some 
$
  x \in \Lambda. 
$
Put 
$       
B_*(\Lambda) = \cup_{k=0}^{\infty} B_k(\Lambda)  
$
where $B_0(\Lambda)$ 
denotes the empty word
$\emptyset$.
Set
$$
\Lambda_l(x) = \{ \mu \in B_l(\Lambda) | \mu x \in X_{\Lambda} \}
\qquad 
\text{ for } x \in X_{\Lambda}, 
\quad
l \in \Zp.
$$
We define a nested sequence of equivalence relations in the space 
$X_{\Lambda}$ as follows (\cite{Ma2}, \cite{Ma4}, \cite{Ma5}).
Two points $x, y \in X_{\Lambda}$ are said to be 
$l$-{\it past equivalent},
written as
$x \sim_l y$,
if $\Lambda_l(x) = \Lambda_l(y)$.
Denote by 
$\Omega_l = X_{\Lambda} / \sim_l$
the quotient space of $X_{\Lambda}$ under $\sim_l$.
For $x,y \in X_{\Lambda}$
 and 
 $\mu \in B_k(\Lambda)$,
 one sees that 
\begin{enumerate}
\renewcommand{\labelenumi}{(\roman{enumi})}
\item if 
$x \sim_l y$, 
we have 
$x \sim_m y$ for $m < l$,
\item if 
$x \sim_l y$ and  $\mu x \in X_{\Lambda}$,
we have 
$ \mu y \in X_{\Lambda}$ and $\mu x \sim_{l-k} \mu y$ for $l>k$.
\end{enumerate}
We have the following sequence of surjections in a natural way:
$$
\Omega_0 \leftarrow
\Omega_1 \leftarrow
\Omega_2 \leftarrow \cdots \leftarrow
\Omega_l \leftarrow
\Omega_{l+1}\leftarrow \cdots,
$$
where $\Omega_0$ is a singleton.
The subshift
 $\Lambda$ is a sofic shift if and only if 
$\Omega_{l} = \Omega_{l+1}$ 
for some $l \in {\Bbb N}$.
For a fixed 
$l \in \Zp$, 
let 
$F_i^l, i= 1, 2,\dots, m(l)$ 
be the set  of all $l$-past equivalence classes of $X_{\Lambda}$
so that $X_{\Lambda}$ is a disjoint union of 
 $F_i^l, i= 1, 2,\dots, m(l)$. 
 Then the canonical $\lambda$-graph system
$
{\frak L}^\Lambda =(V^\Lambda,E^\Lambda,\lambda^\Lambda,\iota^\Lambda )
$
for $\Lambda$ is defined as follows
(\cite{Ma5}). 
 The vertex set $V_l$ at level $l$
consist of the sets 
 $F_i^l,i=1,\dots,m(l)$.
 We write an edge with label $a$ 
    from the vertex 
 $F_i^l \in V_l$ to 
 $F_{j}^{l+1} \in V_{l+1}$ if
$  a x \in F_i^l$
 for some 
 $x \in F_j^{l+1}$.
 We denote by $E_{l,l+1}$ the set of all edges from $V_l$ to $V_{l+1}$.
There exists a natural map $\iota^{\Lambda}_{l,l+1}$ 
from $V_{l+1}$ to
  $V_l$ by mapping $F_j^{l+1}$ to $F_i^l $
 when 
  $F_i^l $ contains  $F_j^{l+1}$.
Set
$V^\Lambda = \cup_{l=0}^{\infty} V_l $
 and
$E^\Lambda = \cup_{l=0}^{\infty} E_{l,l+1}$.
The labeling of edges is denoted by
$\lambda^\Lambda:E \rightarrow \Sigma$.

Let 
$
{\frak L} =(V, E, \lambda, \iota)
$ 
be a $\lambda$-graph system over 
$\Sigma.$ 
The $C^*$-algebra 
${\cal O}_{\frak L}$ 
associated with ${\frak L}$ 
is the  universal    
$C^*$-algebra generated by partial isometries
$S_\alpha, \alpha \in \Sigma$
and projections $E_i^l, i=1,2,\dots,m(l),\ l \in \Zp$  
satisfying  the following  operator relations:
\begin{align}
& \sum_{\alpha \in \Sigma}  S_{\alpha}S_{\alpha}^*  =  1, \\
  \sum_{i=1}^{m(l)} E_i^l  & =  1, \qquad 
 E_i^l   =  \sum_{j=1}^{m(l+1)}I_{l,l+1}(i,j)E_j^{l+1}, \\
& S_\alpha S_\alpha^* E_i^l =E_i^l 
S_\alpha S_\alpha^*,  \\ 
S_{\alpha}^*E_i^l S_{\alpha} & =  
\sum_{j=1}^{m(l+1)} A_{l,l+1}(i,\alpha,j)E_j^{l+1},
\end{align}
for
$
i=1,2,\dots,m(l),\l\in \Zp, 
 \alpha \in \Sigma,
 $
where $V_l = \{v_1^l,\dots, v_{m(l)}^l\},$
\begin{align*}
A_{l,l+1}(i,\alpha,j)
 & =
\begin{cases}
1 &  
    \text{ if } \ s(e) = {v}_i^l, \lambda(e) = \alpha,
                       t(e) = {v}_j^{l+1} 
    \text{ for some }    e \in E_{l,l+1}, \\
0           & \text{ otherwise,}
\end{cases} \\
I_{l,l+1}(i,j)
 & =
\begin{cases}
1 &  
    \text{ if } \ \iota_{l,l+1}({v}_j^{l+1}) = {v}_i^l, \\
0           & \text{ otherwise}
\end{cases} 
\end{align*}
for
$
i=1,2,\dots,m(l),\ j=1,2,\dots,m(l+1), \ \alpha \in \Sigma.
$ 
If $\frak L$ satisfies  $\lambda$-condition $(I)$ 
and is $\lambda$-irreducible,
the $C^*$-algebra ${\cal O}_{\frak L}$
is simple and purely infinite
(\cite{Ma9}, \cite{Ma10}).
By the universality, 
the correspondence $\rho^{\frak L}_z$
$z \in {\Bbb C}$ with $|z| =1$
defined by
$
\rho_z^{\frak L}(S_\alpha) = z S_\alpha,
\rho_z^{\frak L}(E_i^l) = E_i^l
$
for
$\alpha \in \Sigma, i=1,\dots,m(l), l \in \Zp$
yields an action 
$\rho^{\frak L} :
z \in {\Bbb T}
\longrightarrow 
\rho^{\frak L}_z\in\Aut({\cal O}_{\frak L})
$ called gauge action.

Let $\A_{{\frak L},l}$
be the $C^*$-subalgebra of ${\cal O}_{\frak L}$
generated by the projections
$E_i^l, i=1,\dots,m(l)$.
We denote by 
$\A_{\frak L}$ the $C^*$-subalgebra of ${\cal O}_{\frak L}$
generated by the all projections
$E_i^l, i=1,\dots,m(l), l\in \Zp.$
Let 
$\lambda_{\frak L}$ 
be the positive operator on the algebra $\A_{\frak L}$
defined  by
$$
\lambda_{\frak L}(X) = \sum_{\alpha \in \Sigma} S_\alpha^* X S_\alpha \qquad 
\text{ for }
\quad
X \in \A_{\frak L}.
$$ 
Let 
$\lambda_l:\A_{{\frak L},l} \longrightarrow \A_{{\frak L},l+1}$
be the restriction of
$\lambda_{\frak L}$ 
to
$\A_{{\frak L},l}$
and
$\iota_l:\A_{{\frak L},l} \hookrightarrow \A_{{\frak L},l+1}$
the natural inclusion.

We will in this paper study the $C^*$-algebras 
${\cal O}_{{\frak L}^{Z_N}}$
associated with the canonical $\lambda$-graph systems
${\frak L}^{Z_N}$
for nonsofic subshifts $Z_N, N\in {\Bbb N}$
defined in the next section.
The $C^*$-algebras 
${\cal O}_{{\frak L}^{Z_N}}$
will be denoted by
$\OZN$.

\section{The $C^*$-algebra  $\OZN$}
For a finite set $\Sigma$, 
denote by $\Sigma^*$ the set of all words of $\Sigma$.
A (finite or infinite) collection ${\cal C} (\subset \Sigma^*)$
of words over $\Sigma$ is said to be uniquely decipherable 
if
whenever
$\alpha_1 \alpha_2 \cdots \alpha_n 
=  
\gamma_1 \gamma_2 \cdots \gamma_m$
with
$\alpha_i, \gamma_j \in {\cal C}$,
then $n=m$ and $\alpha_i = \gamma_i$ for $i=1,\dots,n$.  
A uniquely decipherable set ${\cal C}$ 
is called a code.
Blanchard and Hansel \cite{BH}
have introduced the notion of coded system.
A subshift $\Lambda$ is called a {\it coded system}
  if $\Lambda$ is the closure of the set of biinfinite sequences obtained by freely 
  concatenating the words in a code ${\cal C}$.
It is denoted by $\Lambda_{\cal C}$.
In this section, 
we will study certain coded systems written as $Z_N, N \in {\Bbb N}$.
We fix a natural number $N \in {\Bbb N}$.
Let $\Sigma_N$ be a set 
$\{ c,b,\alpha_1,\dots,\alpha_N  \}$
of symbols.
The subshift $Z_N$ is defined to be the subshift
 over $\Sigma_N$ whose forbidden words are
$$
 {\cal F}_N = \{ \alpha_i b^m c^k \alpha_j \mid i,j=1,\dots, N,\ m,k =0,1,\dots \text{ with } m \ne k  \}
 $$
 where the word 
 $ \alpha_i b^m c^k \alpha_j $
 means 
  $
  \alpha_i\underbrace{b\cdots b}_{\text{$m$ times}}
  \underbrace{c \cdots c}_{\text{$k$ times}} \alpha_j.  
 $
 It is not a sofic subshift
 and hence not a Markov shift.
Put 
\begin{equation*}
W_{b,c} = \{  w \in \{ b, c \}^* \mid w \ne b^m c^k \text{ for } m,k =0,1,\dots \text{ with } m \ne k  \}
\end{equation*}
the set of all words of $b, c$ that are not of the form 
$b^m c^k $ for $ m,k =0,1,\dots $
with 
$ m \ne k.$ 
\begin{prop}
The set
$
{\cal C}_N := \{ \alpha_i w \mid  i=1,\dots, N, \ w \in W_{b,c}\}
$
 is a code such that $\Lambda_{{\cal C}_N} = Z_N$. 
\end{prop} 
\begin{pf}
It is clear that 
${\cal C}_N $ is a code.
The inclusion relation 
$\Lambda_{{\cal C}_N} \subset Z_N$
is obvious. 
Conversely, the forbidden words
$ {\cal F}_N$  are not admissible in 
$\Lambda_{{\cal C}_N}$. 
Hence
 $\Lambda_{{\cal C}_N} = Z_N$.
 \end{pf}
Define sequences of subsets of $X_{Z_N}$
as in the following way.
\begin{align*}
P_0 = & \{ c^k b^{\infty} \mid k \ge 0 \} \cup \{ 
         b^k c^m b y \in X_{Z_N} \mid k \ge 0, m \ge 1, y \in X_{Z_N} \}\\
\intertext{and}
  E_n = & \{ c^n \alpha_j y \in X_{Z_N} \mid y \in X_{Z_N},\ j=1,\dots,N \}, \\
  Q_l = & \cup_{n > l}E_n, \\
  F_n = & \{ b^m c^{m+n} \alpha_j y \in X_{Z_N} \mid m \ge 1, y \in X_{Z_N},\ j=1,\dots,N \},\\
  R_l = & \{ b^m c^k \alpha_j y \in X_{Z_N} \mid m \ge 1, k \ge 0, m + n \ne k 
            \text{ for  } n = 0,1,\dots,l,  j=1,\dots,N \}
\end{align*}
 for $ l,n = 0,1,\dots.$
\begin{lem}
For each $l \in \Bbb N$,
 the space $X_{Z_N}$ 
is decomposed 
as a disjoint union: 
$$
X_{Z_N} = P_0  \sqcup_{n=0}^{l-1} E_n \sqcup Q_{l-1}  
           \sqcup_{n=0}^{l-1} F_n \sqcup R_{l-1}.
$$
This decomposition of $X_{Z_N}$ 
into $2l+3$-components  
corresponds to the $l$-past equivalence classes of 
$X_{Z_N}$.
\end{lem}
Irreducibility and condition (I) for square matrices with entries in $\{0,1\}$
have been generalized to
$\lambda$-graph systems as 
$\lambda$-irreducibility and $\lambda$-condition (I)
respectively (\cite{Ma10}).
\begin{lem}
The canonical $\lambda$-graph system ${\frak L}^{Z_N}$ 
satisfies $\lambda$-condition (I) and $\lambda$-irreducible.
\end{lem}
\begin{pf}
For any natural number $l$, 
set
 $
L = l + 4.
$ 
For a word $\nu \in B_{l+2}(Z_N),$
put $\mu = \nu c b \in B_L(Z_N)$. 
Then we have for $x \in X_{Z_N}$
\begin{align*}
 \mu x  \in E_k \quad & \text{ 
if } \quad \nu = c^k \alpha_j^{l-k+2} \ \text{ for some } j=1,\dots,N 
         \quad ( k=0,1,\dots,l-1), \\
 \mu x  \in Q_{l-1} \quad & \text{ if } \quad \nu = c^l \alpha_j^2 \ \text{ for some } j=1,\dots,N, \\ 
 \mu x  \in F_k \quad & \text{ if } \quad \nu = b c^{k+1} \alpha_j^{l-k} \ \text{ for some } j=1,\dots,N 
         \quad ( k=0,1,\dots,l-1), \\
 \mu x  \in R_{l-1} \quad & \text{ if } \quad \nu = b^2 c \alpha_j^{l-1} \ \text{ for some } j=1,\dots,N, \\ 
 \mu x  \in P_0 \quad & \text{ if } \quad \nu = c b^2.
\end{align*}
Since the family
$
E_0,F_0,
E_1,F_1,\dots
E_{l-1},F_{l-1},
Q_{l-1}, R_{l-1},P_0
$ 
represents the set of all $l$-past equivalence classes,
the canonical $\lambda$-graph system ${\frak L}^{Z_N}$
for the subshift $Z_N$ 
is $\lambda$-irreducible.
It is easy to see that  ${\frak L}^{Z_N}$ satisfies $\lambda$-condition (I).
\end{pf} 
Therefore we conclude by Lemma 3.3,

\begin{cor}
The $C^*$-algebra $\OZN$ is simple and purely infinite.
\end{cor}

The $C^*$-algebra $\OZN$ is generated by $N+2$ partial isometries
$S_c, S_b, S_{\alpha_1},\dots,S_{\alpha_N}$.
We set 
$T_1 := S_c, T_2 := S_b, S_j: = S_{\alpha_j}, j=1,\dots,N$.
Since for any $x \in X_{Z_N}$, 
both $c x$ and $b x$ are admissible and hence 
belong to
$X_{Z_N}$,
both the operators  
$T_1$ and $T_2$
 are isometries.
Since
for any $x \in X_{Z_N}$,
$\alpha_i x \in X_{Z_N}$ if and only if
$\alpha_j x \in X_{Z_N}$ for $i,j = 1,\dots,N$,
one has 
$S_i^* S_i = S_j^* S_j$ for $i, j=1,\dots,N$.
It has been proved in 
\cite{Ma6} 
that for a subshift $\Lambda$ in general, the associated $C^*$-algebra 
${\cal O}_{\Lambda}$
 can be realized as a universal $C^*$-algebra 
as in the following way:

\begin{lem}[\cite{Ma6}, cf. \cite{CaMa}]
For a subshift $\Lambda$ over $\Sigma = \{1,2,\dots,n \}$,
the $C^*$-algebra ${\cal O}_{\Lambda}$ 
associated with $\Lambda$
 is the universal concrete $C^*$-algebra generated by 
 $n$ partial isometries $S_i, i= 1,2,\dots,n$
subject to the following relations:
\begin{enumerate}\renewcommand{\labelenumi}{(\roman{enumi})}
\item
$
\sum_{j=1}^n S_jS_j^* =1,
$
\item
$
S_i^*S_i = 1 - \sum_{k=1}^{\infty} \sum_{\nu \in L_i^k}
           S_{\nu}S_{\nu}^*, \qquad i=1,2,\dots,n
$
\end{enumerate}
where
 $ L_i^k = \{ \nu_1 \cdots \nu_k \in B_k(\Lambda) \mid
                i\nu_1 \cdots \nu_{k-1} \in B_k(\Lambda), 
                i\nu_1 \cdots \nu_{k-1} \nu_k \not\in
                B_{k+1}(\Lambda) \}.
$  
The infinite sum of the right hand side of the relation \text{(ii)} is taken
under strong operator topology on a Hilbert space.
\end{lem}
The above lemma means that there exists a representation of 
${\cal O}_{\Lambda}$ in operators on a Hilbert space such that the canonical generators satisfy the
relations both (i) and (ii).
Conversely, 
if there exist $n$ partial isometries on a Hilbert space satisfying the above relations, 
then there exists a canonical surjective homomorphism from 
${\cal O}_{\Lambda}$ to the $C^*$-algebra generated by them.  

We apply the preceding lemma to our $C^*$-algebra $\OZN$.
The following lemma is clear.  
\begin{lem}
$
L_{ \alpha_i  }^k =
 \{ b^m c^l \alpha_j \mid  m+l = k-1, m\ne l, \ \text{and } j=1,\dots,N \}
 $
 for $i=1,\dots,N$ and
 $
 k=2,3,\dots.
$
\end{lem}
Thus we obtain 
\begin{thm}
For $N\in {\Bbb N}$,
the $C^*$-algebra $\OZN$ 
associated with the subshift $Z_N$ is simple and purely infinite.
It is the universal concrete 
$C^*$-algebra generated by two isometries $T_1, T_2$
and $N$ partial isometries $S_j, j=1,\dots,N$ 
satisfying the following relations:
\begin{enumerate}\renewcommand{\labelenumi}{(\roman{enumi})}
\item
$
T_1 T_1^* + T_2 T_2^* + \sum_{j=1}^N S_jS_j^* =1,
$
\item
$
S_i^*S_i = 1 - \sum_{k=1}^{\infty} \sum_{m=0, k \ne 2m}^{k} \sum_{j=1}^N
           T_{2^m1^{k-m}} S_j S_j^* T_{2^m1^{k-m}}^*,\qquad i=1,\dots,N
$
\end{enumerate}
where
 $ T_{2^m1^{k-m}}$ 
 denotes
  $
  \underbrace{T_2 \cdots T_2}_{\text{$m$ times}}
  \underbrace{T_1 \cdots T_1}_{\text{$k-m$ times}}.  
 $
The infinite sum of the right hand side of the relation \text{(ii)} is taken
under strong operator topology on a Hilbert space.
\end{thm}

\section{The K-theory for $\OZN$}
In this section, we compute the K-groups for $\OZN$.
Let $(M_{l,l+1},I_{l,l+1})_{l\in \Zp}$
be the nonnegative matrix system
for the symbolic matrix system 
$(\M_{l,l+1}, I_{l,l+1})_{l\in \Zp}$
of the canonical $\lambda$-graph system
${\frak L}^{Z_N}$.
Let $m(l)$ be the row size of the matrix $M_{l,l+1}$.
The main tool is the following K-theory formulae proved in \cite{Ma2}, \cite{Ma9}.
The formulae hold for the 
$C^*$-algebras associated with $\lambda$-graph systems in general.
\begin{lem}{(\cite{Ma2}, \cite{Ma9})}
\begin{enumerate}\renewcommand{\labelenumi}{(\roman{enumi})}
\item
$
K_0(\OZN ) 
=  
\varinjlim 
\{ {\Bbb Z}^{m(l+1)} / (M_{l,l+1}^t - I_{l,l+1}^t){\Bbb Z}^{m(l)}, \overline{I^t}_{l,l+1} \}. 
$
\item
$
K_1(\OZN ) 
=  
\varinjlim \{ \Ker(M_{l,l+1}^t - I_{l,l+1}^t)
             \text{ in }  {\Bbb Z}^{m(l)}, I^t_{l,l+1} \}. 
$
\end{enumerate}
where $M_{l,l+1}^t, I_{l,l+1}^t $ 
are the transposes of the matrices
$M_{l,l+1}, I_{l,l+1} $
respectively,
and
$\overline{I^t}_{l,l+1}:{\Bbb Z}^{m(l)} / (M_{l-1,l}^t - I_{l-1,l}^t){\Bbb Z}^{m(l-1)}
\longrightarrow 
{\Bbb Z}^{m(l+1)} / (M_{l,l+1}^t - I_{l,l+1}^t){\Bbb Z}^{m(l)}
$
is the natural quotient map induced by
$I_{l,l+1}^t$. 
\end{lem}
We will first find  
the symbolic matrix system 
$(\M_{l,l+1}, I_{l,l+1})_{l\in \Zp}$
of the canonical $\lambda$-graph system
${\frak L}^{Z_N}$.
We have
\begin{enumerate}\renewcommand{\labelenumi}{(\arabic{enumi})}
\item $\alpha_i P_0 \subset P_0$ and $ \alpha_i E_0, \alpha_i F_0 \subset E_0$ for $i=1,\dots,N$,
\item $b P_0 \subset P_0, \ \  b E_0, b F_0, b Q_l, b R_l \subset R_{l-1}$
and $b E_n, b F_n  \subset F_{n-1}$ for $n=1,2,\dots, l$,
\item $c P_0, c F_n, c R_l \subset P_0$ for
$n=0,1,\dots,l,\ \ c E_n \subset E_{n+1}$ for $n=0,1,\dots,l-2$
and
$c E_{l-1}, c E_l, c Q_l \subset Q_{l-1}$,   
\end{enumerate}
and
$$
 Q_{l-1} = E_l \sqcup Q_l \quad \text{ and } \quad R_{l-1} = F_l \sqcup R_l.
$$
The we can represent 
the transposes 
$\M_{l,l+1}^t$ 
and
$I_{l,l+1}^t$ 
of the matrices
$\M_{l,l+1}$ 
and
$I_{l,l+1}$ 
respectively as 
$$
\setcounter{MaxMatrixCols}{16}
\M_{l,l+1}^t
=
\begin{bmatrix}
b+c    & \alpha_1+\cdots + \alpha_N & 0 &  \hdotsfor{10}                                       &0& 0 \\
0      & \alpha_1+\cdots + \alpha_N & 0 &c & 0 &  \hdotsfor{8}                                 &0& b \\
c      & \alpha_1+\cdots + \alpha_N & 0 &  \hdotsfor{10}                                       &0& b \\
0      & 0 & b &0 & 0 & c & 0 & \hdotsfor{7}                          & 0 \\
c      & 0 & b &0 & \hdotsfor{10}                                     & 0 \\
0      & 0 & 0 &0 & b & 0 & 0 & c & 0     & \hdotsfor{5}              & 0 \\
c      & 0 & 0 &0 & b & 0 & \hdotsfor{8}                              & 0 \\
0      & 0 & 0 &0 & 0 & 0 & b & 0 & 0     &  c   &   0   &\hdotsfor{3}& 0 \\
c      & 0 & 0 &0 & 0 & 0 & b & 0 &       \hdotsfor{6}                & 0 \\
\vdots &   &   &  &   &   & \ddots  &\ddots   &\ddots &      &       &\ddots  & & & \vdots  \\
\vdots &   &   &  &   &   &   &   &\ddots &\ddots &\ddots &          &\ddots    & & \vdots  \\
0      &   \hdotsfor{8}               &  0   &   b   &   0        &0&c& 0 \\
c      & 0 &   \hdotsfor{7}           &  0   &   b   &   0        &0&0& 0 \\
0      &   \hdotsfor{10}                             &   0        &b&c& 0 \\
c      & 0 & \hdotsfor{9}                            &   0        &b&0& 0 \\
0      & 0 & \hdotsfor{10}                                        &0&c&b  \\
c      & 0 &   \hdotsfor{10}                                      &0&0&b
\end{bmatrix}, 
$$
and
$$
I_{l,l+1}^t 
=
\begin{bmatrix}
1 &0 &   \hdotsfor{12}                                          &0  \\
0 &1 & 0 &   \hdotsfor{10}                              &\dots  &0  \\
0 &0 & 1 & 0 &   \hdotsfor{10}                                  &0  \\
0 &0 & 0 & 1 & 0 &   \hdotsfor{9}                               &0  \\
0 &0 & 0 & 0 & 1 & 0 &   \hdotsfor{8}                           &0  \\
0 &0 & 0 & 0 & 0 & 1 & 0 &   \hdotsfor{7}                       &0  \\
0 &0 & 0 & 0 & 0 & 0 & 1 &0  &   \hdotsfor{6}                   &0  \\
0 &0 & 0 & 0 & 0 & 0 & 0 &1  & 0 &    \hdotsfor{5}              &0  \\
0 &0 & 0 & 0 & 0 & 0 & 0 &0  & 1 & 0  &   \hdotsfor{4}          &0  \\
0 &0 & 0 & 0 & 0 & 0 & 0 &0  & 0 & 1  & 0 &       \hdotsfor{3}  &0  \\
\vdots &  &   &   &   &   &   &   &   & \ddots &\ddots   &\ddots &   &   & \vdots  \\
\vdots &  &   &   &   &   &   &   &   &    &\ddots&  \ddots & \ddots&       & \vdots  \\
\vdots &  &   &   &   &   &   &   &   &    &   &  \ddots &\ddots &\ddots & \vdots  \\
0 &  \hdotsfor{10}                        & \dots & 0     & 1   &0  \\
0 &  \hdotsfor{12}                                        & 0   &1  \\
0 &  \hdotsfor{10}                        & \dots & 0     & 1   &0 \\
0 &  \hdotsfor{12}                                        & 0   &1
\end{bmatrix} 
$$
along the  ordered basis
\begin{align*}
(P_0,E_0,F_0,E_1,F_1,\dots,E_{l}, F_{l}, Q_{l},R_{l})&-\text{row}\\
(P_0,E_0,F_0,E_1,F_1,\dots,E_{l-1}, F_{l-1}, Q_{l-1},R_{l-1})&-\text{column}.  
\end{align*}

The commutative $C^*$-algebras 
$\A_{{\frak L}^{Z_N},l}$
and $\A_{{\frak L}^{Z_N}}$
are denoted by
$\A_{Z_N,l}$ 
and
$\A_{Z_N}$
respectively.
The set of the minimal projections
$\A_{Z_N,l}$ 
correspond to
the set of all $l$-past equivalence classes
$
P_0,
E_0,F_0,E_1, F_1,
\dots,
E_{l-1},F_{l-1},
Q_{l-1}, 
R_{l-1}
$ 
which are denoted by
$
p_0,
e_0,f_0,e_1,f_1,
\dots,
e_{l-1},f_{l-1},
q_{l-1}, r_{l-1}
$ 
respectively.
Hence
$ \text{dim}(\A_l) = m(l) = 2 l +3$  for $l\ge 2$.
We can then represent the induced  matrix on the $K_0$-groups 
$$
M_{l,l+1}^t (=\lambda_{l*}) : 
K_0(\A_{Z_N,l}) = {\Bbb Z}^{2l+3} \rightarrow K_0(\A_{Z_N,l+1}) = {\Bbb Z}^{2l+5}
$$
of the operator $\lambda_l$ 
to be  the matrix $M_{l,l+1}^t$
as follows:
$$
\setcounter{MaxMatrixCols}{16}
M_{l,l+1}^t
=
\begin{bmatrix}
2      & N & 0 &  \hdotsfor{10}                                       &0& 0 \\
0      & N & 0 &1 & 0 &  \hdotsfor{8}                                 &0& 1 \\
1      & N & 0 &  \hdotsfor{10}                                       &0& 1 \\
0      & 0 & 1 &0 & 0 & 1 & 0 & \hdotsfor{7}                          & 0 \\
1      & 0 & 1 &0 & \hdotsfor{10}                                     & 0 \\
0      & 0 & 0 &0 & 1 & 0 & 0 & 1 & 0     & \hdotsfor{5}              & 0 \\
1      & 0 & 0 &0 & 1 & 0 & \hdotsfor{8}                              & 0 \\
0      & 0 & 0 &0 & 0 & 0 & 1 & 0 & 0     &  1   &   0   &\hdotsfor{3}& 0 \\
1      & 0 & 0 &0 & 0 & 0 & 1 & 0 &       \hdotsfor{6}                & 0 \\
\vdots &   &   &  &   &   & \ddots  &\ddots   &\ddots &      &       &\ddots  & & & \vdots  \\
\vdots &   &   &  &   &   &   &   &\ddots &\ddots &\ddots &          &\ddots    & & \vdots  \\
0      &   \hdotsfor{8}               &  0   &   1   &   0        &0&1& 0 \\
1      & 0 &   \hdotsfor{7}           &  0   &   1   &   0        &0&0& 0 \\
0      &   \hdotsfor{10}                             &   0        &1&1& 0 \\
1      & 0 & \hdotsfor{9}                            &   0        &1&0& 0 \\
0      & 0 & \hdotsfor{10}                                        &0&1&1  \\
1      & 0 &   \hdotsfor{10}                                      &0&0&1
\end{bmatrix}. 
$$
The natural inclusion
$\iota_l$
from $\A_{Z_N,l}$ to $\A_{Z_N,l+1}$ 
induces the matrix 
$$
I_{l,l+1}^t (=\iota_{l*}) : 
K_0(\A_{Z_N,l}) = {\Bbb Z}^{2l+3} \rightarrow K_0(\A_{Z_N,l+1}) = {\Bbb Z}^{2l+5}.
$$
As
$$
\setcounter{MaxMatrixCols}{16}
M_{l,l+1}^t -I_{l,l+1}^t
=
\begin{bmatrix}
1      & N & 0 &  \hdotsfor{10}                                       &0& 0 \\
0      &N-1& 0 &1 & 0 &  \hdotsfor{8}                                 &0& 1 \\
1      & N &-1 &0 & \hdotsfor{9}                                      &0& 1 \\
0      & 0 & 1 &-1& 0 & 1 & 0 & \hdotsfor{7}                          & 0 \\
1      & 0 & 1 &0 &-1 & 0 & \hdotsfor{8}                              & 0 \\
0      & 0 & 0 &0 & 1 &-1 & 0 & 1 & 0     & \hdotsfor{5}              & 0 \\
1      & 0 & 0 &0 & 1 & 0 &-1 & 0 & \hdotsfor{6}                      & 0 \\
0      & 0 & 0 &0 & 0 & 0 & 1 & -1& 0     &  1   &   0   &\hdotsfor{3}& 0 \\
1      & 0 & 0 &0 & 0 & 0 & 1 & 0 & -1 & 0& \hdotsfor{4}              & 0 \\
\vdots &   &   &  &   &   & \ddots  &\ddots   &\ddots &      &       &\ddots  & & & \vdots  \\
\vdots &   &   &  &   &   &   &   &\ddots &\ddots &\ddots &          &\ddots    & & \vdots  \\
0      &   \hdotsfor{8}               &  0   &   1   &   -1      &0 &1& 0 \\
1      & 0 &   \hdotsfor{7}           &  0   &   1   &   0       &-1&0& 0 \\
0      &   \hdotsfor{10}                             &   0       &1 &0& 0 \\
1      & 0 & \hdotsfor{9}                            &   0       &1 &0&-1 \\
0      & 0 & \hdotsfor{10}                                       &0 &0&1  \\
1      & 0 &   \hdotsfor{12}                                          &0
\end{bmatrix}, 
$$
we easily see that
\begin{lem}
$
\Ker(M_{l,l+1}^t - I_{l,l+1}^t) 
= 0
\quad
\text{ for }
\quad 
2 \le l \in \Bbb N.
$
\end{lem}
Thus we have by Lemma 4.1,
\begin{prop}
$
K_1(\OZN )\cong 0. 
$
\end{prop}

We will next compute $K_0(\OZN )$.
For an integer $n$, 
we denote by $q(n) \in {\Bbb Z}$ 
the quotient of $n$ by $N$
and 
by $r(n) \in \{0,1,\dots, N-1\}$
its residue such as
$n = q(n) N + r(n)$.
\begin{lem}
Fix $l=2,3,\dots $.
For
$z =
\begin{bmatrix}
z_1  \\
\vdots   \\
z_{2l+5}
\end{bmatrix} 
\in {\Bbb Z}^{2l+5},
$
put inductively
\begin{align*}
x_{2l+3}  & = z_{2l+4},\\
x_{2l+1}  & = z_{2l+2},\\
x_{2l-1}  & = z_{2l+1}- z_{2l+5} + z_{2l+2},\\
x_{2l-3}  & = z_{2l-1}- z_{2l+5} + x_{2l-1}, \\
x_{2l-5}  & = z_{2l-3}- z_{2l+5} + x_{2l-3}, \\
x_{2l-7}  & = z_{2l-5}- z_{2l+5} + x_{2l-5}, \\
          & \vdots \\
x_3       & = z_5 -z_{2l+5} + x_5, \\
x_1       & =  z_{2l+5}, \\
\intertext{ and }
x_2      &  = q(z_1- z_{2l+5}),\\
x_4      &  =   z_2- z_{2l+4} - (N-1)x_2,\\
x_6      &  = z_4 - x_3 + x_4,\\
x_8      &  = z_6 - x_5 + x_6,\\
         & \vdots \\
x_{2l}   & = z_{2l-2} - x_{2l-3} + x_{2l-2},\\
x_{2l+2} & = z_{2l} - x_{2l-1} + x_{2l}.
\end{align*}   
Set
\begin{align*}
r_{l+1}(z) & = r(z_1 - z_{2l+5}) \in \{0,1,\dots,N-1 \},\\
\eta_{l+1}(z) & = z_3 -(x_1 + Nx_2) + x_3 -  x_{2l+3}\\
                             & = -z_1 + (z_3 + z_5 + \cdots + z_{2l+1}) + z_{2l+2} - z_{2l+4}
                             -(l-1)z_{2l+5} + r_{l+1}(z),\\
\psi_{l+1}(z) & = z_{2l+3} - x_1 - x_{2l+1} + x_{2l+3}\\
                             & = -z_{2l+2} + z_{2l+3} +z_{2l+4} - z_{2l+5}.
\end{align*}
Then we have
\begin{equation*}
\begin{bmatrix}
z_1  \\
\vdots   \\
z_{2l+5}
\end{bmatrix} 
=
(M_{l,l+1}^t - I_{l,l+1}^t)
\begin{bmatrix}
x_1  \\
\vdots   \\
x_{2l+3}
\end{bmatrix} 
+
\begin{bmatrix}
r_{l+1}(z)  \\
0 \\
\eta_{l+1}(z) \\
0 \\
\vdots   \\
0 \\
\psi_{l+1}(z)\\
0\\
0
\end{bmatrix}.
\end{equation*}
\end{lem}
\begin{pf}
For 
$[y_i]_{i=1}^{2l+3}\in {\Bbb Z}^{2l+3}$, 
one sees 
\begin{equation*}
(M_{l,l+1}^t - I_{l,l+1}^t)
\begin{bmatrix}
y_1  \\
\vdots   \\
y_{2l+3}
\end{bmatrix} 
=
\begin{bmatrix}
y_1 + N y_2  \\
(N-1) y_2 + y_4 + y_{2l+3} \\
y_1 + N y_2 - y_3 + y_{2l+3} \\
y_3 -y_4 +y_6 \\
y_1 + y_3 -y_5 \\
\vdots \\
y_{2k-1} -y_{2k} + y_{2k+2}  \\
y_1 + y_{2k-1} -y_{2k+1}  \\
\vdots \\
y_{2l-1} -y_{2l} + y_{2l+2}  \\
y_1 + y_{2l-1} -y_{2l+1} \\
y_{2l+1}\\
y_1 + y_{2l+1} -y_{2l+3}\\
y_{2l+3}\\
y_1\\
\end{bmatrix}.
\end{equation*}
Hence the assertion is straightforward.
\end{pf}
\begin{lem}
For
$z= [z_i]_{i=1}^{2l+5} \in {\Bbb Z}^{2l+5}$,
one has 
$$
r_{l+1}(z) = 0 \text{ in } \{0,1,\dots,N-1 \} 
\quad
\text{ and }
\quad
\eta_{l+1}(z) = \psi_{l+1}(z)=0 \text{ in } {\Bbb Z} 
$$
if and only if
there exists
$y= [y_i]_{i=1}^{2l+3}\in {\Bbb Z}^{2l+3}$ 
such that 
$ z =(M_{l,l+1}^t - I_{l,l+1}^t)(y)$.
\end{lem}
\begin{pf}
The only if part follows from the preceding lemma.
We will show the if part.
Suppose that  
there exists
$y= [y_i]_{i=1}^{2l+3}\in {\Bbb Z}^{2l+3}$ 
such that 
$ z =(M_{l,l+1}^t - I_{l,l+1}^t)(y)$.
Define 
$ [x_i]_{i=1}^{2l+3} \in {\Bbb Z}^{2l+3}$
for $z$ as in the preceding lemma
and put
$$
\rho_{l+1}(z) =
[r_{l+1}(z),0,\eta_{l+1}(z),0,\dots,0,\psi_{l+1}(z),0,0]^t
\in {\Bbb Z}^{2l+5}.
$$
Put
$w_i = y_i - x_i, 1\le i \le 2l+3$
so that 
the equality
\begin{equation}
(M_{l,l+1}^t - I_{l,l+1}^t)([w_i]_{i=1}^{2l+3}) = \rho_{l+1}(z)
\end{equation}
holds.
From the $2l+5$-th row in the equality (4.1), 
one sees that
$w_1 =0$.
From the $2l+4$-th row and the $2l+2$-th row, one sees that
$w_{2l+3} =w_{2l+1} =0$
so that $\psi_{l+1}(z) =0$.
Inductively, from the $2k+1$-row, one sees that
$w_{2k-1} = w_{2k+1} -w_1 =0$ for $k=l,l-1, \dots,2$.
From the first row with $w_1 =0$,
one has
$N w_2 = r_{l+1}(z) \in \{0,1,\dots,N-1\}$
so that 
$w_2 =0$ and $r_{l+1}(z) =0$.  
Hence one has from the third row
$\eta_{l+1}(z) = w_1 + N w_2 -w_3 + w_{2l+3} =0$.
\end{pf}
\begin{lem}
The map
$
\xi_{l+1} : 
[z_i]_{i=1}^{2l+5} \in {\Bbb Z}^{2l+5} \longrightarrow
(r_{l+1}(z), \eta_{l+1}(z), \psi_{l+1}(z)) \in \{0,1,\dots,N-1 \}\oplus{\Bbb Z}\oplus{\Bbb Z} 
$
induces an isomorphism from
$
{\Bbb Z}^{2l+5}/(M_{l,l+1}^t - I_{l,l+1}^t){\Bbb Z}^{2l+3}
$
onto
$ 
{\Bbb Z}/N{\Bbb Z} \oplus {\Bbb Z}\oplus{\Bbb Z}.
$
\end{lem}
\begin{pf}
It suffices to show the surjectivity of $\xi_{l+1}$.
For 
$(g,m,k) \in \{0,1,\dots,N-1 \}\oplus{\Bbb Z}\oplus{\Bbb Z}$,
put
$z =[g,0,m,0,\dots,0,k,0,0]^t \in {\Bbb Z}^{2l+5}$.
One then sees that
$$
r_{l+1}(z) = g, \qquad \eta_{l+1}(z) =m, \qquad  \psi_{l+1}(z) =k.
$$
\end{pf}
We denote by 
$\bar{\xi}_{l+1}$
the above isomorphism
from
$
{\Bbb Z}^{2l+5}/(M_{l,l+1}^t - I_{l,l+1}^t){\Bbb Z}^{2l+3}
$
onto
$ 
{\Bbb Z}/N{\Bbb Z} \oplus {\Bbb Z}\oplus{\Bbb Z}
$
induced by 
${\xi}_{l+1}$.
\begin{lem}
The diagram 
$$
\begin{CD}
{\Bbb Z}^{2l+3}/(M_{l-1,l}^t - I_{l-1,l}^t){\Bbb Z}^{2l+1}
@>\overline{I^t}_{l,l+1}>> 
{\Bbb Z}^{2l+5}/(M_{l,l+1}^t - I_{l,l+1}^t){\Bbb Z}^{2l+3}
 \\
@V \bar{\xi}_l VV @V \bar{\xi}_{l+1} VV \\
{\Bbb Z}/N{\Bbb Z} \oplus {\Bbb Z}\oplus{\Bbb Z} 
@>L >> 
{\Bbb Z}/N{\Bbb Z} \oplus {\Bbb Z}\oplus{\Bbb Z}
\end{CD}
$$
is commutative,
where
$L
=
\begin{bmatrix}
1 & 0 & 0 \\
0 & 1 & 1 \\
0 & 0 & 0 
\end{bmatrix}.
$
\end{lem}
\begin{pf}
For 
$(g,m,k) \in \{0,1,\dots,N-1 \}\oplus{\Bbb Z}\oplus{\Bbb Z}$,
put 
$z =[g,0,m,0,\dots,0,k,0,0]^t \in {\Bbb Z}^{2l+3}$.
Since
$
I^t_{l,l+1}(z) =[g,0,m,0,\dots,0,k,0,0,0,0]^t \in {\Bbb Z}^{2l+5},
$
one sees that
$$
r_{l+1}(I^t_{l,l+1}(z)) = g, \qquad \eta_{l+1}(I^t_{l,l+1}(z)) =m+k, \qquad  \psi_{l+1}(I^t_{l,l+1}(z)) =0.
$$
Hence the above diagram is commutative.
\end{pf}

Therefore we  conclude
\begin{thm}
 $K_0(\OZN ) \cong {\Bbb Z}/N{\Bbb Z} \oplus {\Bbb Z}$,
and the class $[1]$
of the unit $1$ of the algebra $\OZN$ in  $K_0(\OZN )$
is $ 0 \oplus 0$ in $ {\Bbb Z}/N{\Bbb Z} \oplus {\Bbb Z}$.
\end{thm}
\begin{pf}
By Lemma 4.1, 
it follows that
\begin{align*}
K_0(\OZN ) 
= & \varinjlim 
\{ {\Bbb Z}^{2l + 5} / (M_{l,l+1}^t - I_{l,l+1}^t){\Bbb Z}^{2l+3}, \overline{I^t}_{l,l+1} \} \\
= & \varinjlim \{ {\Bbb Z}/N{\Bbb Z} \oplus {\Bbb Z}\oplus{\Bbb Z}, L \} \\
\cong & {\Bbb Z}/N{\Bbb Z} \oplus {\Bbb Z}.
\end{align*}
The class of the unit $[1]$ in $K_0(\OZN)$ 
corresponds to 
$[1]$ in $K_0(\A_{{\frak L},l})$
that is represented to be the vector
$[1,\dots,1] \in {\Bbb Z}^{2l+3}$.
Since
$r_{l+1}([1]) = \eta_{l+1}([1])= \psi_{l+1}([1]) = 0$,
we see
$\xi_{l+1}([1]) = (0,0,0).$
Thus the projection $[1]$ in $K_0(\OZN )$ 
represents  $0\oplus 0$ in 
${\Bbb Z}/N{\Bbb Z} \oplus {\Bbb Z} $.
\end{pf}

By \cite{C3}, we know
$K_0({\cal O}_{\infty}) = {\Bbb Z},
K_1({\cal O}_{\infty}) = 0$
and
the position of 
$[1]$ 
in 
$K_0({\cal O}_{\infty})$
is $1$ in ${\Bbb Z}$.
The classification theorem of purely infinite simple $C^*$-algebras proved by
Kirchberg \cite{Kir} and Philips \cite{Ph},
we have 

\begin{cor}
\begin{enumerate}\renewcommand{\labelenumi}{(\roman{enumi})}
\item
The $C^*$-algebra
$\OZN$ for $N\ge 2$ is not stably isomorphic to any Cuntz-Krieger algebra nor to
the  Cuntz algebra ${\cal O}_{\infty}$ 
of order infinity.  
\item
The $C^*$-algebra
$\OZN$ for $N=1$ is not  isomorphic to any Cuntz-Krieger algebra nor to  
the Cuntz algebra ${\cal O}_{\infty}$ 
of order infinity but stably isomorphic to ${\cal O}_{\infty}$.
In fact, ${\cal O}_{Z_1}$ is isomorphic to
$({\cal O}_{\infty})_{1-{s_1s_1^*}}$
where
${\cal O}_{\infty}$ is generated by a sequence 
$s_i, i=1,2,\dots$
of isometries with mutually orthogonal ranges.
\end{enumerate}
\end{cor}

\section{KMS-states for gauge action and topological entropy}
In this section, 
we study KMS-state for gauge action on $\OZN$
and compute topological entropy for the subshift $Z_N$.

In  studying KMS-state for gauge action on the 
$C^*$-algebras associated with subshifts,
the lemma bellow proved in \cite{MWY} plays a crucial r\^{o}le 
(cf. \cite{LaNe}, \cite{PWY}).

\begin{lem}[\cite{MWY}]
For a $\lambda$-graph system ${\frak L}$, if a state 
$\varphi$ on $\OL $ 
is a KMS-state for gauge action on $\OL$ 
at the inverse temperature $ \log \beta$
for some $1< \beta \in {\Bbb R}$,
it satisfies the condition
$\varphi \circ \lambda_{\frak L} = \beta \varphi$
on 
$\A_{\frak L}$.
Conversely,
a state
$\varphi$
on $\A_{\frak L}$ satisfying the condition
$\varphi \circ \lambda_{\frak L} = \beta \varphi$
can be uniquely extended to a KMS-state for gauge action on $\OL$.
\end{lem} 

We will apply this lemma to our $C^*$-algebra $\OZN$ 
for the subshift $Z_N$. 
Let us find a state on $\A_{Z_N}$ satisfying the condition
$\varphi \circ \lambda_{Z_N} = \beta \varphi$
for some 
real number
$\beta $ with $1<\beta \in {\Bbb R}$,
where $\lambda_{Z_N}$
denotes 
$\lambda_{{\frak L}^{Z_N}}$.

Let
$
\hat{p}_0, 
\hat{e}_0, 
\hat{f}_0,  
\hat{e}_1, 
\hat{f}_1,  
\dots,
\hat{e}_l, 
\hat{f}_l,  
\hat{q}_l,\hat{r}_l
$
be  real numbers for $l \in {\Bbb N}$
satisfying the condition
\begin{equation}
 \qquad \hat{p}_0 + 
 \sum_{j=0}^{l}(\hat{e}_j + \hat{f}_j) + \hat{q}_l + \hat{r}_l = 1. 
\end{equation}
We consider the following
equations :
$$
M_{l,l+1}^t
\begin{bmatrix}
\hat{p}_0\\
\hat{e}_0\\
\hat{f}_0\\
\hat{e}_1\\
\hat{f}_1\\
\hat{e}_2\\
\hat{f}_2\\
\vdots\\
\hat{e}_{l-1}\\
\hat{f}_{l-1}\\
\hat{e}_l\\
\hat{f}_l\\
\hat{q}_l\\
\hat{r}_l
\end{bmatrix}
=
\beta
\begin{bmatrix}
\hat{p}_0\\
\hat{e}_0\\
\hat{f}_0\\
\hat{e}_1\\
\hat{f}_1\\
\hat{e}_2\\
\hat{f}_2\\
\vdots\\
\hat{e}_{l-1}\\
\hat{f}_{l-1}\\
\hat{q}_{l-1}\\
\hat{r}_{l-1}
\end{bmatrix}
\quad
\text{ and }
\quad
I_{l,l+1}^t
\begin{bmatrix}
\hat{p}_0\\
\hat{e}_0\\
\hat{f}_0\\
\hat{e}_1\\
\hat{f}_1\\
\hat{e}_2\\
\hat{f}_2\\
\vdots\\
\hat{e}_{l-1}\\
\hat{f}_{l-1}\\
\hat{e}_l\\
\hat{f}_l\\
\hat{q}_l\\
\hat{r}_l
\end{bmatrix}
=
\begin{bmatrix}
\hat{p}_0\\
\hat{e}_0\\
\hat{f}_0\\
\hat{e}_1\\
\hat{f}_1\\
\hat{e}_2\\
\hat{f}_2\\
\vdots\\
\hat{e}_{l-1}\\
\hat{f}_{l-1}\\
\hat{q}_{l-1}\\
\hat{r}_{l-1}
\end{bmatrix} 
$$
for a real number $\beta > 1$.
This means that
\begin{align}
2 \hat{p}_0 + \sum_{i=0}^{l} \hat{f}_i + \hat{r}_l 
& = \beta \hat{p}_0,\\
N (\hat{p}_0 + \hat{e}_0 + \hat{f}_0) 
& = \beta \hat{e}_0,\\
 \hat{e}_{n} + \hat{f}_{n} 
& = \beta \hat{f}_{n-1}, \quad n =1,2, \dots,l, \\
\hat{e}_{n-1} 
& = \beta \hat{e}_{n}, \quad n =1,2, \dots,l-1, \\
\hat{e}_{l-1} + \hat{e}_{l} + \hat{q}_{l} 
& = \beta \hat{q}_{l-1}, \\
\hat{e}_{0} + \hat{f}_{0} +\hat{q}_{l} + \hat{r}_{l}
& = \beta \hat{r}_{l-1} 
\end{align}
and the equations 
\begin{align}
\hat{e}_l + \hat{q}_l & = \hat{q}_{l-1}, \\
\hat{f}_l + \hat{r}_l & = \hat{r}_{l-1}. 
\end{align}
These equations can be solved 
as in the following way by straightforward calculation.
Put
\begin{align*}
g_N(\beta) & = (\beta -1 )\beta \{ (2N-1) \beta -(N+1) \},\\
h_N(\beta) & = (\beta +1 ) 
\{ (\beta -N) (\beta -2) (\beta-1)^2 \beta + (N-2)(\beta -1 )\beta
- N( \beta - 1) -1  \}\\
\intertext{ and }
F_N(\beta) & = h_N(\beta) - g_N(\beta).
\end{align*}
\begin{lem}
\begin{align}
\hat{p}_0 & = \frac{1}{N} \{(\beta - N) \hat{e}_0 - N \hat{f}_0   \},\\ 
(\beta -1)& (\beta+1) g_N(\beta) \hat{f}_0  = h_N(\beta) \hat{e}_0,\\
\hat{e}_n & = \frac{\beta -2}{ \beta^{n+1}},\\
\hat{f}_n & = 
\beta^n \hat{f}_0 - 
(\beta^{n-1}\hat{e}_1 + \beta ^{n-2} \hat{e}_2 + \cdots + \beta \hat{e}_{n-1} 
+ \hat{e}_n ), \\
\hat{q}_n  & = \frac{\beta -2}{ \beta^{n+1}(\beta - 1)}, \\ 
\hat{r}_n  & = 
(\beta -2)\hat{p}_0 - (\hat{f}_0 + \hat{f}_1 + \cdots + \hat{f}_n),
\quad n=0,1,2\dots l. 
\end{align}
\end{lem}
We henceforth assume that the real numbers
$
\hat{p}_0, 
\hat{e}_0, 
\hat{f}_0,  
\hat{e}_1, 
\hat{f}_1,  
\dots,
\hat{e}_l, 
\hat{f}_l,  
\hat{q}_l,\hat{r}_l
$
are all nonnegative for all
$l \in {\Bbb N}$.
\begin{lem}
\begin{enumerate}\renewcommand{\labelenumi}{(\roman{enumi})}
\item $ \beta >2$ and $\hat{e}_n, \hat{q}_n >0$ for $n=0,1,\dots,l$. 
\item  $\beta > N$.  
\item $g_N(\beta) >0$.
\end{enumerate}
\end{lem}
\begin{pf}
(i) 
Since $\hat{e}_1\ge 0$, one sees that 
$\beta \ge 2$ by (5.12).
Suppose that $\beta =2$.
One has $\hat{e}_n = \hat{q}_n = 0$ for $n=0,1,\dots,l$.
(5.15) implies that $\hat{r}_n = \hat{f}_n = 0$ for $n=0,1,\dots,l$.
It follows 
that $\hat{p}_0 = 0$ by (5.10), a contradiction to (5.1).
Hence $\beta >2$ and $\hat{e}_n, \hat{q}_n >0$ for $n=0,1,\dots,l$.

(ii)
Suppose next that $\beta \le N$.
(5.10) implies that $\hat{p}_0 = 0,  \hat{f}_0 =0$.
Then (5.13) implies that $\hat{f}_n <0$ a contradiction.

(iii)
As
\begin{equation}
\frac{g_N(\beta)}{(\beta-1)\beta} = (2N-1)\beta -(N+1) 
> (2N-1)N - (N+1) =(N-1)^2 + N^2 -2.
\end{equation}
one sees that
$g_N(\beta) >0$ for $N\ge 2$ and
$g_1(\beta) = (\beta -2)(\beta-1) \beta >0$
\end{pf}
\begin{lem}
The conditions
$\hat{f}_n \ge 0$ for all $n = 0,1,\dots$
imply
$
F_N(\beta) \ge 0.
$
\end{lem}
\begin{pf}
(5.12) and (5.13) imply 
\begin{equation}
\hat{f}_n  = 
\beta^n 
\{ \hat{f}_0 - \frac{1}{\beta^2 -1}( 1- \frac{1}{\beta^{2n}} ) \hat{e}_0 \}
\end{equation}
so that the conditions
$
\hat{f}_n \ge 0 
$
for all 
$
n =0, 1, \dots
$
imply         
\begin{equation}
(\beta^2-1) \hat{f}_0 -\hat{e}_0 \ge 
 - \frac{1}{\beta^{2n}}  \hat{e}_0 \quad \text{ for all } n =0, 1, \dots
\end{equation}
As $g_N(\beta) >0$, 
by (5.11), the above equalities imply the the inequality
$
h_N(\beta) - g_N(\beta) \ge 0.
$ 
\end{pf}
\begin{lem}
The conditions 
$\hat{r}_n \ge 0$ for all $n = 0,1,\dots.$
imply $F_N(\beta)\le 0$.
\end{lem}
\begin{pf}
By (5.12) and (5.13),
one sees that
\begin{equation*}
\sum_{k=0}^n \hat{f}_k 
= \frac{\beta^{n+1}-1}{\beta-1} \hat{f}_0 
+ 
\{ \frac{1}{ (\beta -1)^2 } 
- \frac{1}{(\beta -1)^2 (\beta+1)}(\beta^{n+1} + \frac{1}{\beta^n} )
\}
\hat{e}_0.
\end{equation*}
(5.15) implies that
\begin{align*}
& (\beta -1)^2 (\beta+1) g_N(\beta)\hat{r}_n \\
=
& [\frac{\beta-N}{N}(\beta -2)(\beta -1)^2 (\beta +1) 
  - \{ (\beta +1) - (\beta^{n+1} + \frac{1}{\beta^n}) \} ] 
  g_N(\beta) \hat{e}_0 \\
&  -[(\beta -2)(\beta-1)^2 (\beta+1) + (\beta-1)(\beta +1) (\beta^{n+1}-1)] 
g_N(\beta) \hat{f}_0 \\
=
& [\frac{\beta-N}{N}(\beta -2)(\beta -1)^2 (\beta +1) 
  - \{ (\beta +1) - (\beta^{n+1} + \frac{1}{\beta^n}) \} ] 
  g_N(\beta) \hat{e}_0 \\
&  -[(\beta -2)(\beta-1) +  (\beta^{n+1}-1)] h_N(\beta) \hat{e}_0 \\
= 
& - F_N(\beta) \beta^{n+1} \hat{e}_0\\
& + [\frac{\beta-N}{N}(\beta -2)(\beta -1)^2 (\beta +1) 
  - \{ (\beta +1) -  \frac{1}{\beta^n} \} ] 
  g_N(\beta) \hat{e}_0 \\
&  -[(\beta -2)(\beta-1) -1 ] h_N(\beta) \hat{e}_0 
\end{align*}
Since
$$
\liminf_{n\to{\infty}}\frac{
(\beta -1)^2 (\beta+1) g_N(\beta)\hat{r}_n }
{\beta^{n+1}}
=
\liminf_{n\to{\infty}}\frac{
(\beta -1)^2 (\beta+1) g_N(\beta)\hat{r}_n }
{\beta^{n+1}}
= - F_N(\beta) \hat{e}_0
\ge 0
$$
and
$\hat{e}_0 >0$,
we see
$F_N(\beta)\le 0$.
\end{pf}
Therefore we have 
\begin{cor}
Suppose that
$
\hat{p}_0,
\hat{e}_n,
\hat{f}_n,
\hat{q}_n, 
\hat{r}_n 
$
are all nonnegative for all $n = 0,1,\dots$.
Then $\beta >N$ and 
$F_N(\beta) =0$.
\end{cor}

Conversely, for $N \ge 2$, we have
\begin{lem}
For $\beta >N$, 
$F_N(\beta) =0$ implies 
$\hat{f}_n, \hat{r}_n >0$ for all $n=0,1,\dots $.
\end{lem}
\begin{pf}
Assume that 
$F_N(\beta) =0$.
As $\beta >N$,
(5.16) implies 
$g_N(\beta) >0$.
Since 
$h_N(\beta) = g_N(\beta)$,
the inequality (5.18) hold,
which means 
$\hat{f}_n > 0$
for all 
$n=0,1,\dots$ by (5.17).

We will next show that
$\hat{r}_n >0$ for all $n=0,1,\dots $.
By the equality 
 $g_N(\beta) = h_N(\beta)$, 
one has
\begin{align*}
 & (\beta -1)^2 (\beta+1) g_N(\beta)\hat{r}_n \\ 
=
 & [\frac{\beta-N}{N}(\beta -2)(\beta -1)^2 (\beta +1) 
  - \{ (\beta +1) -  \frac{1}{\beta^n} \} ] 
  g_N(\beta) \hat{e}_0 \\
&  -[(\beta -2)(\beta-1) -1 ] g_N(\beta) \hat{e}_0. 
\end{align*}
It follows that
\begin{align*}
&  [\frac{\beta-N}{N}(\beta -2)(\beta -1)^2 (\beta +1) 
  - (\beta +1) ] 
  -[(\beta -2)(\beta-1) -1 ]\\
=
& \frac{1}{N\beta}(\beta^3 -N +1) >0
\end{align*}
so that 
$(\beta -1)^2 (\beta+1) g_N(\beta)\hat{r}_n$
is positive.
As $g_N(\beta) >0$,
one sees $\hat{r}_n >0$ for all $n=0,1,\dots $.
\end{pf}
Therefore we have
\begin{cor}
For $\beta >N$, 
$F_N(\beta) =0$ if and only if  
$
\hat{p}_0,
\hat{e}_n,
\hat{f}_n,
\hat{q}_n, 
\hat{r}_n 
$
are all nonnegative for all $n = 0,1,\dots,$.
In this case,
they are all positive.
\end{cor}
\begin{pf}
For $N=1$, the statement has been shown in 
\cite{Ma3}.
\end{pf}
We note that 
the identity 
\begin{equation*}
F_N(\beta) = \beta^6 -(N+3)\beta^5 +(3N+1)\beta^4- 2(N-1)\beta^3 -(N+2)\beta^2 + N-1
\end{equation*}
holds.
For $N=1$, the above identity goes to
\begin{equation*}
F_1(\beta) = \beta^6 - 4\beta^5 + 4\beta^4 - 3\beta^2. 
\end{equation*}
The unique positive solution of 
$ F_1(\beta) = 0$
is 
$\beta = 1 + \sqrt{ 1 + \sqrt{3}} = 2.652891\cdots$
( \cite{Ma3}).

We will next study positive solutions of the equations 
\begin{equation*}
F_N(\beta) =0 \qquad \text{ for } \quad N \ge 2.
\end{equation*}
\begin{lem}
\begin{enumerate}\renewcommand{\labelenumi}{(\roman{enumi})} \qquad \qquad \qquad
\item $F_2(3) < 0 < F_2(4)$.
\item $F_N(N) < 0 < F_N(N+1)$ for $N\ge 3$.
\end{enumerate}
\end{lem}
\begin{pf}
By the identities
\begin{align*}
F_N(N) & = -(N^2 +1)(N^2 -N +1), \\ 
F_N(N+1)& = \{ N(N-2)(N+1) -1 \}(N+1)^2 + N-1
\end{align*}
one sees that (iii) holds.
\end{pf}
The following lemma is straightforward.
\begin{lem}
\begin{enumerate}\renewcommand{\labelenumi}{(\roman{enumi})} \qquad \qquad \qquad
\item For $N=2$,
$F^\prime_2(\beta) >0$ for all $\beta \ge 3$.
Hence 
$F_2(\beta) =0$ has a unique solution 
in $[3,\infty)$.
\item  For $N=3$,
$F^\prime_3(3) <0$ and $F^{\prime \prime}_3(\beta) >0$ for all $\beta \ge 3$.
Hence 
$F_3(\beta) =0$ has a unique solution 
in $[3,\infty)$.
\item  For $N \ge 4$,
$F^\prime_N(\beta) >0$ for all $\beta \ge N$.
Hence 
$F_N(\beta) =0$ has a unique solution 
in $[N,\infty)$.
\end{enumerate}
\end{lem}
Therefore we have
\begin{prop}
For $N \ge 1$,
the equation $F_N(\beta) =0$
has a unique solution in $[N,\infty)$.
It satisfies $N < \beta < N+1$ for $N \ge 3$.  
\end{prop}
\begin{lem}
For the  unique $\beta >N$ satisfying 
$F_N(\beta) =0$, one has 
\begin{align*}
\hat{p}_0 & = \frac{(\beta -2)(\beta^2-N\beta -1)}{N(\beta^2 - 1)},\\ 
\hat{e}_n & = \frac{\beta -2}{ \beta^{n+1}},\\
\hat{f}_n & = \frac{\beta -2}{(\beta^2 - 1) \beta^{n+1}},\\
\hat{q}_n  & = \frac{\beta -2}{ \beta^{n+1}(\beta - 1)}, \\ 
\hat{r}_n  & = \frac{\beta-2}{N(\beta-1)^2(\beta+1)}\{(\beta-2)(\beta-1)(\beta^2- N\beta -1) -N + \frac{1}{\beta^{n+1}}\}. 
\end{align*}
\end{lem}
We remark that 
for $N=1$, under the condition 
$F_1(\beta) =0$, one sees
$
\hat{p}_0 = \frac{1}{(\beta - 1)^2}.
$

Let
$
p_0,
e_0, f_0,  
e_1, f_1, \dots,
e_{l-1}, f_{l-1},  
q_{l-1}, r_{l-1}
$
be 
the set of minimal projections of $\A_{Z_N,l}$ considered in the previous section.
For a state $\varphi$ on $\A_{Z_N}$,
put
$$
\hat{t}:=\varphi(t) 
\qquad
\text{for}
\quad
t=
p_0, \
e_0, \ f_0, \  
e_1, \ f_1, \ \dots,
e_{l-1}, \ f_{l-1}, \  
q_{l-1}, \ r_{l-1}.
$$
Then
$\varphi$ satisfies the condition
$\varphi \circ \lambda_{Z_N} = \beta \varphi$
 on $\A_{Z_N}$ 
if and only if the real numbers
$
\hat{p}_0, 
\hat{e}_0, \hat{f}_0,  
\hat{e}_1, \hat{f}_1,\dots,
\hat{e}_{l-1}, \hat{f}_{l-1},  
\hat{q}_{l-1}, \hat{r}_{l-1}
$
are all nonnegative and satisfy all the equations from 
(5.1) to (5.9)
for all $l = 2,3,\dots$.
Therefore we have
\begin{prop}
A state $\varphi$ on  $\A_{Z_N}$ 
satisfies the condition
$\varphi \circ \lambda_{Z_N} = \beta \varphi$
for some 
real number
$\beta $
 if and only if
 $\beta >N$
and $\beta$ is the unique solution of the equation:
\begin{equation*}
\beta^6 -(N+3)\beta^5 +(3N+1)\beta^4- 2(N-1)\beta^3 -(N+2)\beta^2 + N-1=0.
\end{equation*}
Moreover
a state $\varphi$ that satisfies the above condition 
is unique.
\end{prop}

\def\topen{{ h_{\text{top}}(\Lambda) }}

The topological entropy
of some classes of subshifts,
including irreducible topological Markov shifts,
  $\beta$-shifts for real number $\beta > 1$
and the subshift $Z_N$ for $N=1$
 have appeared 
as the logarithm of the inverse temperature of the admitted KMS states for gauge actions on 
the associated $C^*$-algebras (\cite{EFW2}, \cite{KMW}, \cite{Ma3}).
For a subshift $(\Lambda, \sigma)$ 
and a natural number $k$,
let $\theta_k(\Lambda)$ be the cardinal number of the words 
$B_k(\Lambda)$ of length $k$
appearing in $\Lambda$.
The topological entropy 
 $\topen$ for $(\Lambda, \sigma)$ 
is given by
$$
  \topen = \lim_{k \to \infty} \frac{1}{k}\text{ log } \theta_k(\Lambda)
  \qquad(cf.\cite{LM}).
$$
 
 For the subshifts $Z_N$, 
we have
\begin{lem}
If there exists a $\log \beta$-KMS state on $\OZN$ for gauge action
for some real number $\beta > N$, 
we have
$$
       \log \beta = \log  r(\lambda_{Z\N})
                        = h_{\text{top}}(Z_N)
$$
where
$r(\lambda_{Z_N})$
 denotes the spectral radius of the operator 
$\lambda_{Z_N}$ on $\A_{Z_N}$.
\end{lem}
\begin{pf}
The proof is similar to the proof of Lemma 6.7 of \cite{Ma3}.
For the sake of completeness, 
we will give a proof.
A word $\mu = \mu_1\dots \mu_n$ in $\Sigma_N$ 
appears in the subshift $Z_N$ if and only if 
$S_{\mu}(= S_{\mu_1}\cdots S_{\mu_n})  \ne 0$.
Let 
$\varphi$ be a log $\beta$-KMS state  on $\OZN$
 for gauge action for some positive real number
$\beta$.
For $k \in \Bbb N$,
it follows that 
$$
\beta^k = \varphi(\sum_{\mu \in B_k(Z_N)}S_{\mu}^*S_{\mu})
        \le \| \sum_{\mu \in B_k(Z_N)}S_{\mu}^*S_{\mu} \|
        = \| \lambda_{Z_N}^k(1) \|
        \le  \sum_{\mu \in B_k(Z_N)} \| S_{\mu}^*S_{\mu} \|
        = \theta_k(Z_N ).
$$
As $\lambda_{Z_N}^k$ 
is a completely positive map on the unital 
$C^*$-algebra $\A_{Z_N}$, 
we have
$\| \lambda_{Z_N}^k(1) \| = \| \lambda_{Z_N}^k \|$
so that we see
\begin{equation}
\beta^k \le \| \lambda_{Z_N}^k \| \le \theta_k({Z_N}).
\end{equation}
On the other hand, 
by  the inequality 
$
     \beta^k 
 \ge \theta_k(Z_N) \min_{\mu \in Z^k} \varphi(S_\mu^* S_\mu),
$
we obtain
$$
       \min_{\mu \in B_k(Z_N)} \varphi(S_\mu^* S_\mu)^{\frac{1}{k}} \cdot
       \theta_k(Z_N)^{\frac{1}{k}} 
  \le
       \beta
  \le  
       \theta_k(Z_N)^{\frac{1}{k}}. 
$$
Now we have
$S_\mu^* S_\mu \ge P_0$
for any word
$\mu \in B_*(Z_N)$.
It follows that
$$
\varphi(S_\mu^* S_\mu) \ge \varphi(P_0)= \frac{(\beta-2)(\beta^2 - N \beta -1)}{N (\beta^2 -1)}>0
\quad
\text{ for }
\quad
\mu \in B_*(Z_N).
$$
Hence we obtain
$$
\lim_{k \to \infty} 
\min_{\mu \in B_k(Z_N)}\varphi(S_\mu^* S_\mu)^{\frac{1}{k}} 
 = 1
$$ 
and 
$
\lim_{k \to \infty} 
       \theta_k(Z_N)^{\frac{1}{k}}
= \beta.
$
Thus we get the desired equalities from (5.19).
\end{pf}
Therefore we conclude
\begin{thm}
\begin{enumerate}\renewcommand{\labelenumi}{(\roman{enumi})}
\item
For a positive real number $\beta$,
there exists a $ \log \beta$ KMS-state for gauge action on  $\OZN$ 
if and only if 
$\beta >N$ and 
$\beta$ is the unique solution of the equation:
\begin{equation*}
\beta^6 -(N+3)\beta^5 +(3N+1)\beta^4- 2(N-1)\beta^3 -(N+2)\beta^2 + N-1=0.
\end{equation*}
\item
The above KMS-state is unique.
\item
$\log \beta =
h_{\text top}(Z_N)$ :
the topological entropy for the subshift $Z_N$.
\end{enumerate}
\end{thm} 

\medskip
We finally mention an asymptotic behavier of the solution
$\beta > N$ for $F_N(\beta) =0.$
Let $\beta_N >N$ 
be the unique solution for the equation
$F_N(\beta) =0$.
We know them by numerical calculation
as in the following way:
\begin{align*}
\beta_1 & = 2.652891650\dots \\
\beta_2 & = 3.063607825\dots \\
\beta_3 & = 3.670666991\dots \\
\beta_4 & = 4.446202651\dots \\
\beta_5 & = 5.321226229\dots \\
\beta_6 & = 6.247124025\dots \\
\beta_7 & = 7.199582119\dots \\
\beta_8 & = 8.166942400\dots \\
        & \dots 
\end{align*}
\begin{prop}
Let
$\beta_N >N$ be the unique solution of the equation
$F_N(\beta) =0$.
Then we have
\begin{equation*}
\lim_{N \to \infty} \frac{\beta_N}{N} =1.
\end{equation*} 
\end{prop}
\begin{pf}
Put
$t_N =\frac{\beta_N}{N}$.
It follows that
\begin{align*}
0 & = \frac{F_N(\beta_N)}{N^6}\\
  & =
 {t_N}^6 -(1 + \frac{3}{N}) {t_N}^5 +(\frac{3}{N}+\frac{1}{N}) {t_N}^4
- 2( \frac{1}{N^2} - \frac{1}{N^3}) {t_N}^3 
 -( \frac{1}{N^3} + \frac{2}{N^4}) {t_N}^2 + \frac{1}{N^5} - \frac{1}{N^6}
\end{align*}
and hence 
\begin{equation*}
1 - t_N  
= \frac{3}{N} + (\frac{3}{N}+\frac{1}{N})\frac{1}{ t_N}
- 2( \frac{1}{N^2} - \frac{1}{N^3}) \frac{1}{{t_N}^2} 
-( \frac{1}{N^3} + \frac{2}{N^4}) \frac{1}{{t_N}^3}
+( \frac{1}{N^5} - \frac{1}{N^6})\frac{1}{{t_N}^5}.
\end{equation*}
As
$
0 < \frac{1}{t_N} <1, 
$
we get $\lim_{N\to{\infty}} | t_N -1 | =0.$
\end{pf}

\section{Flow equivalence classes of the subshifts $Z_N, N \in {\Bbb N}$}
We will finally  apply our discussions to 
a classification problem in symbolic dynamical systems under flow equvalence
(cf. \cite{BF}, \cite{Fr}, \cite{PS}).
In \cite{Ma7}, \cite{Ma8},
we have defined the K-groups 
$K_i(\Lambda), i=0,1$
and the Bowen-Franks groups
$BF^i(\Lambda), i=0,1$
for subshift $\Lambda$
by the K-groups
$K_i({\cal O}_\Lambda), i=0,1$
and the Ext-groups
$\Ext^i({\cal O}_\Lambda), i=0,1$
for the associated $C^*$-algebra
${\cal O}_\Lambda$
respectively.
We have then proved that the groups
$K_i(\Lambda), BF^i(\Lambda), i=0,1$
are invariant under not only topological conjugacy class 
but also flow equivalence class of subshifts.
Especially there is no known computable invariant under flow equivalence of subshifts
other than the groups.
Since the Ext-groups 
$\Ext^i({\cal O}_\Lambda), i=0,1$
can be computed by the Universal Coefficient Theorem for the K-theory 
of $C^*$-algebra ${\cal O}_{\Lambda}$,
we have by Theorem 4.8,
\begin{align*}
K_0(Z_N) = {\Bbb Z}/N{\Bbb Z} \oplus {\Bbb Z},& \qquad K_1(Z_N) = \{ 0\},\\
BF^0(Z_N) = {\Bbb Z}/N{\Bbb Z},& \qquad BF^1(Z_N) = {\Bbb Z}.
\end{align*}

Therefore  we have
\begin{prop}
The subshifts $Z_N, N=1,2,\dots$
are not flow equivalent to each other.
\end{prop}

{\it Acknowledgements.}
The author would like to deeply thank Wolfgang Krieger for his many useful conversations.


\begin{thebibliography}{99}




\bibitem{Be}
{\sc M. P. B{\'{e}}al},
{\it Codage Symbolique},
 Masson, Paris
 (1993).




\bibitem{BH}
{\sc F. Blanchard and G. Hansel},
{\it Systems cod{\'e}s
}, Theor.\ Computer Sci.\ 
{\bf 44}
(1986)
pp.\ 17--49.







\bibitem{BF}
{\sc R. Bowen and J. Franks},
{\it Homology for zero-dimensional nonwandering sets}, Ann.\ Math.\ 
{\bf 106}
(1977)
pp.\ 73--92.





\bibitem{CaMa}
{\sc T. M. Carlesen and K. Matsumoto},
{\it Some remarks on the $C^*$-algebras associated with subshifts},
Math. Scand.\
{\bf 95}(2004), 
pp.\ 145-160.


\bibitem{ChS}{\sc  N. ~Chomsky and M. ~P. ~Sch{\"{u}}tzenberger},
{\it The algebraic theory of context-free languages},
Computer programing and formal systems, North-Holland
(1963),
pp.\ 118--161.




\bibitem{C}
{\sc J. Cuntz},
{\it Simple $C^*$-algebras generated by isometries},
 Comm.\ Math.\ Phys.\
{\bf 57}(1977) pp.\ 173--185.



\bibitem{C2}
{\sc J. Cuntz}, 
{\it A class of $C^*$-algebras and topological Markov chains II: reducible chains and the Ext- functor for $C^*$-algebras}, 
Inventions Math.\
{\bf 63}(1980) pp.\ 25--40.


\bibitem{C3}
{\sc J. Cuntz},
{\it K-theory for certain $C^*$-algebras},
 Ann. Math.\
{\bf 117}(1981) pp.\ 181--197.



\bibitem{CK}{\sc J. ~Cuntz and W. ~Krieger},
{\it A class of $C^*$-algebras and topological Markov chains},
 Invent.\ Math.\
 {\bf 56}(1980), pp.\ 251--268.







\bibitem{EFW2}
{\sc M. Enomoto, M. Fujii and Y. Watatani},
{\it KMS states for gauge action on ${\cal O}_A$}, 
Math.\ Japon
{\bf 29} (1984)
pp.\ 607--619.





\bibitem{Fi}
{\sc R. Fischer},
{\it Sofic systems and graphs},
Monats.\ f{\"u}r Math.\
{\bf 80}
(1975)
pp.\ 179--186.


\bibitem{Fr}
{\sc J. Franks},
{\it Flow equivalence of subshifts of finite type}, 
Ergodic Theory Dynam.\ Systems
{\bf 4}
(1984)
pp.\ 53--66.








\bibitem{HU}{\sc J. ~E. ~Hopcroft and J. ~D. ~Ullman},
{\it Introduction to Automata Theory, Languages, and Computation},
Addison-Wesley, Reading                        
(2001).

\bibitem{KMW}
{\sc Y. Katayama, K. Matsumoto and Y. Watatani},
{\it Simple $C^*$-algebras arising from $\beta$-expansion of real numbers},
Ergodic Theory  Dynam.\ Systems
{\bf 18} (1998) pp.\ 937--962.






\bibitem{Kir}
{\sc E. Kirchberg},
{\it The classification of purely infinite  
$C^*$-algebras using Kasparov's theory},
preprint,(1994).





\bibitem{Ki}{\sc B.~P. ~Kitchens},
{\it Symbolic dynamics}, Springer-Verlag, Berlin, Heidelberg and New York
(1998).







\bibitem{Kr4}
{\sc W. ~Krieger},
{\it On sofic systems I},
Israel J. Math.
{\bf 48}(1984), pp.\ 305--330.


\bibitem{Kr5}
{\sc W. ~Krieger},
{\it On sofic systems II},
Israel J. Math.
{\bf 60}(1987), pp.\ 167--176.



\bibitem{Kr6}
{\sc W. ~Krieger},
{\it  On subshifts and topological Markov chains},  
Numbers, information and complexity (Bielefeld 1998),
Kluwer Acad. Publ. Boston MA (2000)
pp.\ 453--472



\bibitem{KM}{\sc W. ~Krieger and K. ~Matsumoto},
{\it Shannon graphs, subshifts and lambda-graph systems},
J.\ Math.\ Soc.\ Japan 
{\bf 54}(2002),
pp.\ 877--900.

%

%
\bibitem{KM4}{\sc W. ~Krieger and K. ~Matsumoto},
{\it  Subshifts from certain one-counter codes}, preprint.







\bibitem{LaNe}
{\sc M. Laca and S. Neshveyev},
{\it KMS-states of quasi-free dynamics on Pimsner algebras}, 
J.\  Funct.\ Anal.\
{\bf 211}
(2004)
pp.\ 457--482.


\bibitem{LM}{\sc D. ~Lind and B. ~Marcus},
{\it An introduction to symbolic dynamics and coding},
 Cambridge University Press, Cambridge
(1995).


\bibitem{Ma}
{\sc K. Matsumoto},
{\it On $C^*$-algebras associated with subshifts},
Internat. J. Math.
{\bf 8}(1997), pp. 357-374.

\bibitem{Ma2}
{\sc K. Matsumoto},
{\it K-theory for  $C^*$-algebras associated with subshifts},
Math. Scand.
{\bf 82}(1998), pp.\ 237-255.



\bibitem{Ma3}
{\sc K. Matsumoto},
{\it A simple  $C^*$-algebra arising from certain subshift},
 J.\ Operator Theory 
{\bf 42}
 (1999),
 pp.\ 351--370.
  
\bibitem{Ma4}
{\sc K. Matsumoto},
{\it Dimension groups for subshifts and simplicity of the associated $C\sp *$-algebras},
 J.\  Math.\ Soc.\ Japan 
 {\bf 51}
  (1999), 
  pp.\ 679--698. 
  
\bibitem{Ma5}
{\sc K. Matsumoto},
{\it Presentations of subshifts and their topological conjugacy invariants}, 
Documenta Math.\ 
{\bf 4}
 (1999),
 pp.\  285--340. 



\bibitem{Ma6}
{\sc K. Matsumoto},
{\it Relations among generators of  $C^*$-algebras associated with subshifts},
Internat. J. Math.
{\bf 10}(1999), pp.\ 385-405.



\bibitem{Ma7}
{\sc K. Matsumoto},
{\it Bowen-Franks groups for subshifts and Ext-groups for $C\sp *$-algebras},
  $K$-Theory 
{\bf 23}
 (2001), 
 pp.\ 67--104.



\bibitem{Ma8}
{\sc K. Matsumoto},
{\it Bowen-Franks groups as an invariant for flow equivalence of subshifts}, 
Ergodic Theory Dynam. Systems 
{\bf 21}
 (2001), 
 pp.\ 1831--1842. 


\bibitem{Ma9}
{\sc K. Matsumoto},
{\it $C\sp *$-algebras associated with presentations of subshifts},
 Documenta Math.\ 
{\bf 7}
 (2002), 
 pp.\ 1--30. 


\bibitem{Ma10}
{\sc K. Matsumoto},
{\it Construction and pure infniteness of $C^*$-algebra associated with 
lambda-graph systems},  Math.\ Scand.\
{\bf 97} (2005),pp.\ 73--89.











\bibitem{MWY}
{\sc K. Matsumoto, Y. Watatani and M. Yoshida},
{\it KMS-states for gauge actions 
       on $C^*$-algebras associated with subshifts}, 
 Math.\ Z.\
{\bf 228}
(1998)
pp.\ 489--509.








\bibitem{PS}
{\sc W. Parry and D. Sullivan},
{\it A topological invariant for flows on one-dimensional spaces}, 
Topology
{\bf 14}
(1975)
pp.\ 297--299.



\bibitem{PWY}
{\sc C. Pinzari, Y. Watatani and K. Yoshida},
{\it KMS-states, entropy and the variational principle in full $C^*$-dynamical systems}, 
Comm.\  Math.\ Phys.\
{\bf 213}
(2000)
pp.\ 231--379.


\bibitem{Ph}
{\sc N. C. Phillips},
{\it A classification theorem for nuclear  purely infinite simple
 $C^*$-algebras},
Doc.\ Math.\
{\bf 5}
(2000)
pp.\ 49--114.







\bibitem{Ro}
{\sc M. R{\o}dom},
{\it Classification of Cunzt-Krieger algebras},
 K-theory {\bf 9}(1995), pp.\  31--58.

\bibitem{Ro2}
{\sc M. R{\o}dom},
{\it Classification of purely infinite simple $C^*$-algebras I}, 
J. Func. Anal.
{\bf 131}(1995), pp.\ 415--458.










\bibitem{We}
{\sc B. Weiss},
{\it Subshifts of finite type and sofic systems}, 
Monats.\ Math.\
{\bf 77},
(1973),
pp.\ 462--474.



\end{thebibliography}
\end{document}